\documentclass[11pt]{article}
\usepackage{enumerate, fullpage}
\usepackage{graphicx}
\usepackage{multirow}
\usepackage{color}
\usepackage{amsmath}
\usepackage{amsfonts}
\usepackage{amssymb}
\usepackage{mathrsfs}
\usepackage{lscape}
\usepackage{xcolor}
\definecolor{lightblue}{RGB}{173,216,230}
\definecolor{powderblue}{RGB}{176,224,230}
\definecolor{liteblue}{RGB}{190,238,244}
\definecolor{liteseagreen}{RGB}{90,242,184}
\definecolor{liteorangered}{RGB}{255,210,190}
\definecolor{liteorangeredplus}{RGB}{255,230,209}
\definecolor{liteseagreen}{RGB}{90,242,184}
\definecolor{darkseagreen}{RGB}{50,121,92}
\usepackage[colorlinks=true,linkcolor=blue,citecolor=darkseagreen]{hyperref}

\newcommand{\Y}{\mathbb{Y}}

\def\norm#1{\|#1\|}

\newcommand{\R}{\mathbb{R}}

\newcommand{\proj}{\mathop{\rm proj}}

\newcommand{\beq}{\begin{equation}}
\newcommand{\eeq}{\end{equation}}

\newcommand{\beqnr}{\begin{eqnarray}}
\newcommand{\eeqnr}{\end{eqnarray}}

\newcommand{\benum}{\begin{enumerate}}
\newcommand{\eenum}{\end{enumerate}}

\newtheorem{DE}{Definition}[section]

\newtheorem{LE}[DE]{Lemma}

\newtheorem{CO}[DE]{Corollary}

\newcommand{\qed}{\mbox{}\hspace*{\fill}\nolinebreak\mbox{$\rule{0.7em}{0.7em}$}}

\advance \topmargin by -\headheight
\advance \topmargin by -\headsep

\evensidemargin \oddsidemargin

\begin{document}
\Large
\begin{center}
{\bf On linear relaxations of OPF problems}\\
Daniel Bienstock and Gonzalo Mu\~{n}oz, Columbia University
\end{center}
\normalsize
\section{Introduction}.
The AC OPF problem is a fundamental software component in the operation of electrical power transmission systems.  For background, see \cite{bergenvittal}.  It can be formulated as a nonconvex, continuous optimization problem.  In routine problem instances, solutions
of excellent quality can be quickly obtained using a variety of methodologies, including
sequential linearization and interior point methods.   Instances involving grids under
stress or extreme conditions can prove significantly more difficult.

The recent work by Lavaei and Low \cite{javad} on semidefinite programming relaxations
has sparked renewed interest in this problem.  See \cite{molzahn} (and references therein) for some cutting-edge approaches.  

\subsection{Our approach} \label{ours}
Here we focus on developing \textit{linear relaxations} to AC OPF problems, in \textit{lifted spaces}, with the primary goal of quickly proving lower bounds and enabling fast, standard
optimization methodologies such as branching and the incorporation of binary variables into
optimization models.  To motivate our approach, let $(P, Q, V^{(2)})$ be a vector that
includes, for each line $km$, the real and reactive power injections 
$P_{km}, P_{mk}, Q_{km}$ and $Q_{mk}$, and for each bus $k$ the squared bus voltage
magnitude $|V_k|^2$, denoted by $V^{(2)}_k$.  Using these variables, 
we first write the OPF problem in 
the following summarized form
\begin{subequations}
\label{opf1}
\begin{eqnarray}
\min && F(P, Q, V^{(2)}) \\
\mbox{Subject to:} && \\
&& A P \, + \, B Q \, + \, C V^{(2)} \ \le \ d \label{simpleconstraints} \\
&& (P, Q, V^{(2)}) \ \in \ \Omega \label{complex}.
\end{eqnarray}
\end{subequations}
Here,
\begin{itemize}
\item In constraints (\ref{simpleconstraints}), $A$, $B$ and $C$ are matrices and $d$ is
a vector, all of appropriate dimension. These constraints describe basic relationships 
such as generator output limits, $(P,Q)$-bus demand statements, and voltage limits.  These are all linear
constraints and thus can be expressed in the form (\ref{simpleconstraints}).
\item Constraints (\ref{complex}) describe the underlying physics, e.g. Ohm's law.  
For example, in 
the rectangular formulation of AC OPF such constraints of course will involve additional
variables (the real and imaginary voltage components at each bus) and bilinear constraints
relating those variables to the vector $(P, Q, V^{(s)})$.
\item In standard OPF problem formulations, the objective $F(P, Q, V^{(2)})$ is typically
the sum of active power generation costs (summed over the generators) a separable
convex quadratic function of the generator outputs.
\end{itemize}
Our basic approach will approximate (\ref{complex}) with linear inequalities obtained
by lifting formulation (\ref{opf1}) to a higher-dimensional space, and running a 
cutting-plane algorithm over that lifted formulation.  By 'lifting' we mean a procedure
that adds new variables (with specific interpretations) and then writes inequalities
that such variables, together with $(P, Q, V^{(s)})$, must satisfy in a feasible solution
to the OPF problem.  
To fix our language, we view the quantities $P_{km}, P_{mk}, Q_{km}, Q_{mk}$ (for each line $km$) and $|V_k|^2$ (for each bus $k$) as \textit{foundational}.  All other variables,
including those that arise naturally from constraint (\ref{complex}) as well as those
that we introduce, will be called \textit{lifted}\footnote{Occasionally we may view
the rectangular voltage coordinates as foundational.}.\\

In the following sections we introduce our lifted variables, as well as the inequalities
that we introduce so as to obtain a convex relaxation of (\ref{complex}).  The inequalities
will be of four types: 
\begin{enumerate}
\item $\Delta$-inequalities, 
\item (active power) loss inequalities, 
\item Circle inequalities
\item Semidefinite cut inequalities.
\end{enumerate}
All these inequalities are convex; some linear and some conic.  In the case of conic
inequalities we rely on outer approximation through tangent cutting planes so as to
ultimately obtain linear formulations as desired.

In Section \ref{tighten} we present a tightening procedure, and in 
Section \ref{mip} we describe the use of linear mixed-integer programming.

\section{Basic inequalities} \label{basic} We consider a line $\{k,m\}$ with series impedance $z = r + j x$ and series admittance 
\begin{eqnarray}
&& y \ \doteq \ z^{-1} = g + j b, \quad \mbox{where} \label{ydef} \\
&& g = \frac{r}{r^2 +x^2}  \quad \mbox{and} \quad    b = -\frac{x}{r^2 +x^2}. \label{gbdef}
\end{eqnarray}
In addition, there will be a shunt admittance  $y^{sh} = g^{sh} + j b^{sh}$, and a transformer with
tap ratio 
\begin{eqnarray}
&& N  \ \doteq \ \tau e^{j \sigma}
\end{eqnarray}
where $\tau$ is the magnitude and $\sigma$ is the phase shift angle. Note that
$r = r(km)$, etc, but for simplicity of 
notation we omit the dependence of line parameters on $km$.

\begin{figure}[htb]
\centering
\begin{center}
\includegraphics[height=1.5in]{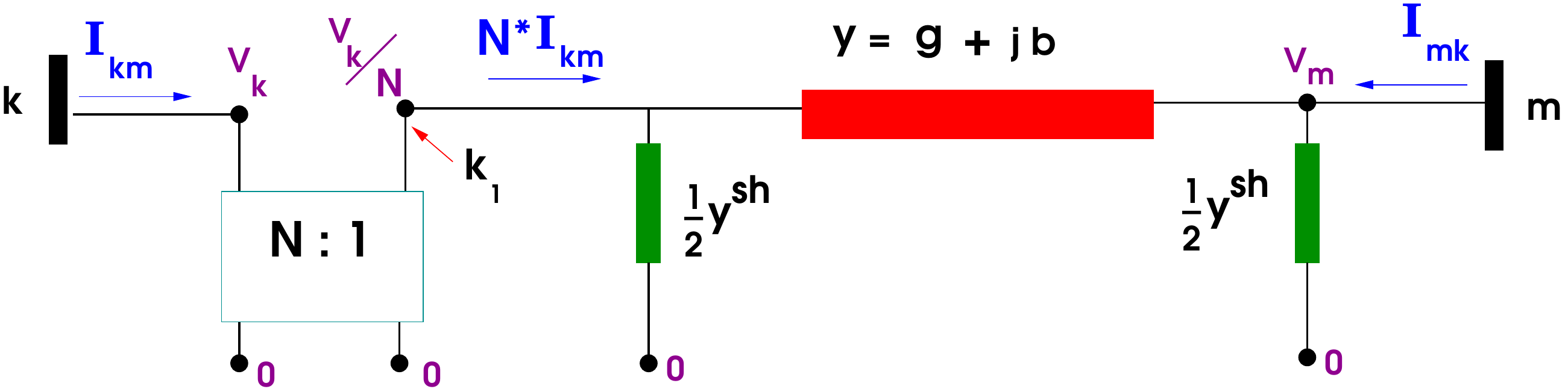}
\vskip -10pt
\caption{$\pi$-model, including transformer and shunt admittance}\label{pimodel}
\end{center}
\end{figure}
\noindent In the figure, voltages are shown in purple and currents in blue.  Notice that the transformer is assumed to be located at the ``k'' or ``from'' end of the line.

\noindent Define:
\begin{eqnarray}
V & = & \left( \begin{array}{r}
V_{k} \\
V_{m}
\end{array} \right) = \left( \begin{array}{r}
|V_k| e^{j \theta_k} \\
|V_m| e^{j \theta_m} \\
\end{array} \right) = \left( \begin{array}{c}
e_k + j f_k \\
e_m + j f_m \\
\end{array} \right) \quad \mbox{(voltages at $k$ and $m$)}\label{voltage}\\
I & = & \left( \begin{array}{r}
I_{km} \\
I_{mk}
\end{array} \right) \quad \quad \mbox{(complex current injections at $k$ and $m$)}\label{current}\\
S & = & \left( \begin{array}{r}
S_{km} \\
S_{mk}
\end{array} \right) \, = \, \left( \begin{array}{r}
P_{km} + j Q_{km}\\
P_{mk} + j Q_{mk}
\end{array} \right) \quad \mbox{(complex power injections at $k$ and $m$)}\label{power}
\end{eqnarray}
Then
\begin{eqnarray}
&& S_{km} \ = \ V_k I^*_{km}, \quad S_{mk} \ = \ V_m I^*_{mk} \quad \mbox{and} \quad I \ = \ \Y V, 
\end{eqnarray}
where $\Y$ is the branch admittance matrix, defined as
\Large
\begin{eqnarray}
\Y & = & \left( \begin{array}{lcr}
(y + \frac{y^{sh}}{2}) \frac{1}{\tau^2} & & -y \frac{1}{\tau e^{-j \sigma}}\\
& \\
-y \frac{1}{\tau e^{j \sigma}} & & y +  \frac{y^{sh}}{2}
\end{array} \right).
\end{eqnarray}
\normalsize
\noindent In the next sections we derive the power equations, the $\Delta$ and the circle inequalities, first for the simplest case (no shunt, no transformer) then for the case with shunts but no transformers, and finally for the most general case.  
\subsection{ $y^{sh} = 0$ and $N = 1$.}
\noindent In this case we have
\begin{eqnarray}
I_{km} & = & y(V_k - V_m). \label{ikm1}
\end{eqnarray}
In rectangular coordinates this means that
\begin{eqnarray}
I_{km} & = & g(e_k - e_m) - b(f_k - f_m)  \, + \, j [\,b(e_k - e_m) + g(f_k - f_m) \, ]
\end{eqnarray}
with a symmetric expression for $I_{mk}$.  Therefore
\begin{eqnarray}
P_{km} & = & e_k g(e_k - e_m) - e_k b(f_k - f_m) + f_k g(f_k - f_m) + f_k b(e_k - e_m) \\
 & = & (e_k - e_m) (g e_k + b f_k) + (f_k - f_m) (-be_k + g f_k)\\
 & = & (e_k - e_m) (g \, , \,  b) (^{e_k}_{f_k}) \ + \ (f_k - f_m) (-b \, , \, g) (^{e_k}_{f_k}) \label{rectpkm1} 
\end{eqnarray}
with a symmetric expression for $P_{mk}$.  Similarly,
\begin{eqnarray}
Q_{km} & = & f_k g(e_k - e_m) - f_k b(f_k - f_m) - e_k g(f_k - f_m) - e_k b(e_k - e_m) \\
 & = & (e_k - e_m) (g f_k - b e_k) + (f_k - f_m) (-g e_k - b f_k). \label{rectqkm1}
\end{eqnarray}
\noindent To obtain similar inequalities in polar coordinates we write the impedance and admittance in polar coordinates:
$$ z = |z|e^{j \angle z}, \quad y = \frac{1}{|z|}e^{-j \angle z}. $$
Then (see e.g. Bergen and Vittal \cite{bergenvittal}, p. 104)
\begin{eqnarray}
S_{km} & = & V_k I^*_{km} \ = \ V_k (V_k^* - V_m^*)y^* \, = \, \frac{|V_k|^2}{|z|}e^{j \angle z} - \frac{|V_k||V_m|}{|z|} e^{j \angle z} e^{j \theta_{km}}, \label{skmpolar1}
\end{eqnarray}
where 
$$ \theta_{km} \doteq \theta_k - \theta_m.$$
We also can rewrite (\ref{skmpolar1}) as
\begin{eqnarray}
S_{km} & = & |V_k|^2(g - j b) - |V_k||V_m|(g - j b)(\cos \theta_{km} + j \sin \theta_{km}) \nonumber \\
& = & |V_k|^2g - |V_k||V_m| g \cos \theta_{km} - |V_k||V_m| b \sin \theta_{km} \nonumber \\
&& + j \left[ -|V_k|^2b + |V_k||V_m| b \cos \theta_{km} - |V_k||V_m| b \sin \theta_{km} \right]. \label{skmpolar1a}
\end{eqnarray}
Likewise, the power received at $m$ (rather than injected), $-S_{mk}$, satisfies
\begin{eqnarray}
-S_{mk} & = & -\frac{|V_m|^2}{|z|}e^{j \angle z} + \frac{|V_k||V_m|}{|z|} e^{j \angle z} e^{-j \theta_{km}}.
\end{eqnarray}
We also obtain an expression for $S_{km}$ similar to (\ref{skmpolar1a}) by
switching the $k$ and $m$ symbols.
\subsubsection {$\Delta$ and loss inequalities, 1}
Let $\mu_{km}$  and $\nu_{km}$ denote  known upper bounds on
$$ | (g \, , \,  b) (^{e_k}_{f_k})| \quad \mbox{and} \quad | (-b \, , \, g) (^{e_k}_{f_k})|,$$
respectively.  Then, using (\ref{rectpkm1}) we obtain:
\begin{eqnarray}
|P_{km}| & \le & \mu_{km} | e_k - e_m|  \, + \, \nu_{km} | f_k - f_m|. \label{delta1}
\end{eqnarray}
This is the basic $\Delta$ inequality.  Note that the vectors $(^g_b)$
and $(^{-b}_g)$ are of equal norm and orthogonal, so further elaborations of 
the $\Delta$-inequalities are possible.

By adding the expression for $P_{km}$ in (\ref{skmpolar1a}) and the
corresponding expression for $P_{mk}$ we obtain
\begin{eqnarray}
P_{km} + P_{mk} & = & g (|V_k|^2 + |V_m|^2) - 2 g |V_k| |V_m| \cos \theta_{km} \ = \ g |V_k - V_m|^2, \label{lossmain}
\end{eqnarray}
\noindent which can be \textit{relaxed} as
\begin{eqnarray}
&& g (e_k - e_m)^2 + g (f_k - f_m)^2  \ \le \ P_{km} + P_{mk}, \label{loss1} 
\end{eqnarray}
or equivalently,
\begin{eqnarray}
&& g |e_k - e_m|^2 + g |f_k - f_m|^2  \ \le \ P_{km} + P_{mk}. \label{loss1b}
\end{eqnarray}
We term (\ref{loss1}) of (\ref{loss1b}) the \textit{loss} inequality.  Note that by definition $g \ge 0$ (unless by a modeling artifact we have $r < 0$).   The point of (\ref{loss1b}) is that in a lifted formulation with variables representing $|e_k - e_m|$ and $|f_k - f_m|$, (\ref{loss1b}) is convex.

\subsubsection {Circle inequalities, 1}
We can rewrite ineq. (\ref{skmpolar1}) as
\begin{eqnarray}
S_{km} & = & C_{km} - B_{km} e^{j \theta_{km}} \quad \mbox{where}\\
C_{km}   & \doteq & \frac{|V_k|^2}{|z|}e^{j \angle z} \ \mbox{and} \\
B_{km}   & \doteq & \frac{|V_k||V_m|}{|z|} e^{j \angle z}.
\end{eqnarray}
\noindent Note that $C_{km}$ and $B_{km}$ are obtained in the 
complex plane by rotating the real 
numbers $\frac{|V_m|^2}{|z|}$ and $\frac{|V_k||V_m|}{|z|}$ (respectively) 
by the same angle $\angle z$.  As $\theta_{km}$ varies, (\ref{skmpolar1})
indicates that $S_{km}$ describes a circle (the ``sending circle'') with center $C_{km}$ and 
radius 
$$\rho \doteq \frac{|V_k||V_m|}{|z|}.$$  
 Likewise, $-S_{mk}$ describes a circle (the ``receiving circle'')
with center $-\frac{|V_m|^2}{|z|}$ and radius $\rho$.
Refer to Bergen and Vittal for more
details.  Using either circle we can obtain valid convex inequalities.   For example, clearly we have
\begin{eqnarray}
&& [ \operatorname{Re}( S_{km} - C_{km} )]^2 \ + \ [ \operatorname{Im}( S_{km} - C_{km} )]^2 \ \le \ \rho^2, \quad \mbox{or in other words} \\
&& \left(P_{km} - \frac{r|V_k|^2}{r^2 + x^2} \right)^2 \ + \ \left(Q_{km} - \frac{x |V_k|^2}{r^2 + x^2} \right)^2 \ \le \ \frac{|V_k|^2 |V_m|^2 }{r^2 + x^2}. \label{firstcircle1}
\end{eqnarray}
As discussed in Section \ref{ours}, our formulation has variables used to represent $P_{km}, Q_{km}, |V_k|^2$ 
and $|V_m|^2$. Using these variables (\ref{firstcircle1}) is a conic constraint.
Using these variables, from (\ref{firstcircle1}) we obtain a convex system by adding
two lifted variables $\alpha_{km}$, $\beta_{km}$ and the constraints
\begin{subequations}
\label{circlesystem1}
\begin{eqnarray}
P_{km} - \frac{r V_k^{(2)}}{r^2 + x^2} & = & \alpha_{km} \\
Q_{km} - \frac{x V_k^{(2)}}{r^2 + x^2} & = & \beta_{km} \\
\alpha_{km}^2 + \beta_{km}^2 & \le & \frac{V_k^{(2)} V_m^{(2)} }{r^2 + x^2}.
\end{eqnarray}
\end{subequations}
\subsection{ General $y^{sh}$ but $N = 1$}
\noindent In this case we have
\begin{eqnarray}
I_{km} & = & y(V_k - V_m) \, + \,  \frac{1}{2} y^{sh} V_k, \label{ikm2}
\end{eqnarray}
and so in rectangular coordinates
\begin{eqnarray}
I_{km} & = & g(e_k - e_m) - b(f_k - f_m)  + \frac{1}{2}(g^{sh}e_k - b^{sh} f_k) \, + \, \nonumber \\
&& \  j [\,b(e_k - e_m) + g(f_k - f_m) + \frac{1}{2} (b^{sh} e_k + g^{sh} f_k)\, ]. \label{ikm2rect}
\end{eqnarray}
We will now obtain 
\begin{eqnarray}
P_{km} & = & (e_k - e_m) (g \, , \,  b) (^{e_k}_{f_k}) \ + \ (f_k - f_m) (-b \, , \, g) (^{e_k}_{f_k}) \, + \, \frac{g^{sh}}{2} (e_k^2 + f_k^2)\label{rectpkm2} 
\end{eqnarray}
and
\begin{eqnarray}
Q_{km} & = (e_k - e_m) (g f_k - b e_k) + (f_k - f_m) (-g e_k - b f_k) - \frac{b^{sh}}{2} (e^2_k + f^2_k). \label{rectqkm2}
\end{eqnarray}
Note that expressions in (\ref{rectpkm2}) and (\ref{rectqkm2}) are obtained
from (\ref{rectpkm1}) and (\ref{rectqkm1}) 
by adding the terms $\frac{g^{sh}}{2} (e^2_k + f^2_k)$ and
$-\frac{b^{sh}}{2} (e^2_k + f^2_k)$, respectively.\\

\noindent To obtain similar expressions under polar coordinates, note that the 
only the expression (\ref{ikm2}) for $I_{km}$ differs from (\ref{ikm1}) only
in the term $\frac{1}{2} y^{sh} V_k$.  Thus, 
\begin{eqnarray}
S_{km} & = & V_k I^*_{km} \ = \ V_k (V_k^* - V_m^*)y^* + \frac{1}{2} (g^{sh} -j b^{sh}) V_k V_k^*. \label{skmpolar2a}
\end{eqnarray}
\subsubsection {$\Delta$ and loss inequalities, 2}
Using (\ref{rectpkm2}), we obtain
\begin{eqnarray}
|P_{km}| - \frac{g^{sh}}{2} V_k^{(2)} & \le & \mu_{km} | e_k - e_m|  \, + \, \nu_{km} | f_k - f_m|. \label{delta2}
\end{eqnarray}
This is the second version of the $\Delta$-inequality\footnote{A common modeling assumption is that $g^{sh} = 0$.}.  
Since the right-hand side of (\ref{rectpkm2}) is obtained 
by adding  $\frac{g^{sh}}{2} (e^2_k + f^2_k)$ to the right-hand side of
(\ref{rectpkm1}), we have the following analogue of (\ref{loss1}):
\begin{eqnarray}
&& g (e_k - e_m)^2 + g (f_k - f_m)^2  \ \le \ P_{km} + P_{mk} - \frac{g^{sh}}{2} (V_k^{(2)} + V_m^{(2)}), \label{loss2}
\end{eqnarray}
the second version of our loss inequality, which again is a conic constraint if we introduce appropriate variables.
\subsubsection {Circle inequalities, 2}
From (\ref{skmpolar2a}) we get
\begin{eqnarray}
S_{km} & = & |V_k|^2\left(\frac{e^{j \angle z} }{|z|} + \frac{1}{2}( g^{sh} - j  b^{sh}) \right) - \frac{|V_k||V_m|}{|z|} e^{j \angle z} e^{j \theta_{km}}, \label{skmpolar2}
\end{eqnarray}
which again describes a circle, with center and radius, respectively,
\begin{eqnarray}
|V_k|^2\left(\frac{e^{j \angle z} }{|z|} + \frac{1}{2}( g^{sh} - j  b^{sh}) \right) \quad \mbox{and} \quad \frac{|V_k||V_m|}{|z|}. \label{centerradius2}
\end{eqnarray}

Using (\ref{skmpolar2}), (\ref{centerradius2}), and since
\begin{eqnarray}
 \operatorname{Re}\left( (|V_k|^2\left(\frac{e^{j \angle z} }{|z|} + \frac{1}{2}( g^{sh} - j  b^{sh}) \right) \right) & = & |V_k|^2 \left( \frac{r }{r^2 + x^2} + \frac{g^{sh}}{2} \right) \ \mbox{and}\\
 \operatorname{Im}\left( (|V_k|^2\left(\frac{e^{j \angle z} }{|z|} + + \frac{1}{2}( g^{sh} - j  b^{sh}) \right) \right) & = &  |V_k|^2 \left( \frac{x}{r^2 + x^2} -  \frac{b^{sh}}{2} \right), 
\end{eqnarray}
we obtain the following generalization of (\ref{firstcircle1}):
\begin{subequations}
\label{circlesystem2}
\begin{eqnarray}
P_{km} - \left( \frac{r }{r^2 + x^2} + \frac{g^{sh}}{2} \right)  V_k^{(2)} & = & \alpha_{km} \\
Q_{km} - \left( \frac{x }{r^2 + x^2}  - \frac{b^{sh}}{2}\right) V_k^{(2)} & = & \beta_{km} \\
\alpha_{km}^2 + \beta_{km}^2 & \le & \frac{V_k^{(2)} V_m^{(2)} }{r^2 + x^2}.
\end{eqnarray}
\end{subequations}

\subsection{ General $b^{sh}$ and $N$}
\noindent In this case we have
\begin{eqnarray}
I_{km} & = & \frac{1}{\tau} y \left[ \frac{1}{\tau} V_k - e^{j \sigma} V_m \right] \, + \,  \frac{1}{2 \tau^2} y^{sh} V_k  \label{ikm3a}\\
& = & \frac{1}{\tau} y \left[ \frac{1}{\tau} V_k - (\cos \sigma + j \sin \sigma) V_m \right] \, + \,  \frac{1}{2 \tau^2} (g^{sh} + j b^{sh}) V_k.  \label{ikm3b}
\end{eqnarray}
\noindent In rectangular coordinates this can be further expanded as
\begin{eqnarray}
&& I_{km} \ = \ \nonumber \\
&& \frac{1}{\tau}(g + jb) \left[ \frac{1}{\tau} (e_k + j f_k) - 
(\cos \sigma + j \sin \sigma) (e_m + j f_m) \right] + 
 \frac{1}{2 \tau^2} (g^{sh} + j b^{sh}) (e_k + j f_k) \nonumber \\
&& = \frac{1}{\tau}(g + jb) \left[ \frac{1}{\tau} (e_k + j f_k) - e_m \cos \sigma + f_m \sin \sigma - j ( e_m \sin \sigma + f_m \cos \sigma) \right] \nonumber \\
&& + \frac{1}{2 \tau^2} \left[ g^{sh}e_k - b^{sh} f_k + j (b^{sh} e_k + g^{sh} f_k) \right] \nonumber \\
&& =  \frac{1}{\tau}(g + jb) \left[ \frac{1}{\tau} e_k - e_m \cos \sigma + f_m \sin \sigma + j ( \frac{1}{\tau} f_k - e_m \sin \sigma  - f_m \cos \sigma) \right] \nonumber \\
&& + \frac{1}{2 \tau^2} \left[ g^{sh}e_k - b^{sh} f_k + j (b^{sh} e_k + g^{sh} f_k) \right].
\end{eqnarray}
From this expression we obtain:
\begin{eqnarray}
\operatorname{Re} I_{km} & = &  \frac{g}{\tau} \left[ \frac{e_k}{\tau} - e_m \cos \sigma + f_m \sin \sigma \right] \nonumber \\
&& - \frac{b}{\tau} \left[ \frac{1}{\tau} f_k - e_m \sin \sigma  - f_m \cos \sigma \right] \nonumber \\
&& + \frac{1}{2 \tau^2} \left[ g^{sh}e_k - b^{sh} f_k \right], \quad \mbox{and} \label{reIkmgen} \\ 
&& \nonumber \\
\operatorname{Im} I_{km} & = &  \frac{b}{\tau} \left[ \frac{e_k}{\tau} - e_m \cos \sigma + f_m \sin \sigma \right] \nonumber \\
&& + \frac{g}{\tau} \left[ \frac{1}{\tau} f_k - e_m \sin \sigma  - f_m \cos \sigma \right] \nonumber \\
&& + \frac{1}{2 \tau^2} \left[ b^{sh} e_k + g^{sh} f_k \right]. \label{imIkmgen}
\end{eqnarray}  
\noindent We can see that in the ``no-transformer'' case, i.e. $\tau = 1$ and $\sigma = 0$, (\ref{reIkmgen}) and (\ref{imIkmgen}) match the expansion
(\ref{ikm2rect}) for $I_{km}$, as desired. We then have:
\begin{eqnarray}
P_{km} & = & \operatorname{Re} V_k I^*_{km} \nonumber \\
  & = & e_k \left\{\frac{g}{\tau} \left[ \frac{e_k}{\tau} - e_m \cos \sigma + f_m \sin \sigma \right] - \frac{b}{\tau} \left[ \frac{1}{\tau} f_k - e_m \sin \sigma  - f_m \cos \sigma \right] \right\} \nonumber \\
&& + f_k \left\{ \frac{b}{\tau} \left[ \frac{e_k}{\tau} - e_m \cos \sigma + f_m \sin \sigma + \right] + \frac{g}{\tau} \left[ \frac{1}{\tau} f_k - e_m \sin \sigma  - f_m \cos \sigma \right]\right\} \nonumber \\
&& + \frac{g^{sh}}{2 \tau^2} (e^2_k + f^2_k). \label{rectpkm3a}
\end{eqnarray}
We can rewrite 
(\ref{rectpkm3a}) as:
\begin{eqnarray}
P_{km} & = & \frac{1}{\tau} \left[ \frac{e_k}{\tau} - e_m \cos \sigma \right] (g \, , \,  b) (^{e_k}_{f_k}) + \frac{1}{\tau} \left[ \frac{1}{\tau} f_k - f_m \cos \sigma \right]  (-b \, , \, g) (^{e_k}_{f_k}) \nonumber \\
&& + \, \frac{g^{sh}}{2 \tau^2} (e^2_k + f^2_k) \nonumber \\
&& + \, \frac{g e_k f_m + b e_k e_m + b f_k f_m - g f_k e_m }{\tau} \sin \sigma . \label{rectpkm3x}
\end{eqnarray}
\noindent In the no-transformer case this expression evaluates to 
\begin{eqnarray}
&& (e_k - e_m) (g \, , \,  b) (^{e_k}_{f_k}) + (f_k - f_m)(-b \, , \, g) (^{e_k}_{f_k}) + \frac{g^{sh}}{2} (e^2_k + f^2_k) 
\end{eqnarray}
which is the same as (\ref{rectpkm2}), as desired.  Note that in
(\ref{rectpkm3x}) the third term vanishes when there is no transformer,
and the second term vanishes when there is no shunt conductance. We can further rewrite (\ref{rectpkm3x}) as
\begin{eqnarray}
P_{km} & = & \frac{1}{\tau} \left[ \frac{e_k}{\tau} - e_m \cos \sigma + f_m \sin \sigma \right] (g \, , \,  b) (^{e_k}_{f_k}) + \frac{1}{\tau} \left[ \frac{1}{\tau} f_k - f_m \cos \sigma - e_m \sin \sigma \right]  (-b \, , \, g) (^{e_k}_{f_k}) \nonumber \\
&& + \, \frac{g^{sh}}{2 \tau^2} (e^2_k + f^2_k). \label{rectpkm3}
\end{eqnarray}

\noindent Next we will compute an expression for $P_{mk}$. In the transformer case the line is not symmetric and we first need to compute $I_{mk}$.  We have:
\begin{eqnarray}
I_{mk} & = & -\frac{1}{\tau e^{j \sigma}} y V_k + \left( y + \frac{y^{sh}}{2} \right) V_m  \label{imk3a} \\
& = & -\frac{1}{\tau}(\cos \sigma - j \sin \sigma)(g + j b)(e_k + j f_k) \nonumber \\ && + \left( g + g^{sh}/2 + j (b + b^{sh}/2) \right) (e_m + j f_m) \nonumber \\
& = & -\frac{1}{\tau}(g + j b) \left[ e_k \cos \sigma + f_k \sin \sigma + j(  - e_k \sin \sigma + f_k \cos \sigma) \right] \nonumber \\
&& + (g + g^{sh}/2) e_m  - (b + b^{sh}/2) f_m  +j \left[ (b + b^{sh}/2) e_m + (g + g^{sh}/2) f_m \right]. \label{imk3b}
\end{eqnarray}
Therefore
\begin{eqnarray}
\operatorname{Re} I_{mk} & = & -\frac{g}{\tau}\left[ e_k \cos \sigma + f_k \sin \sigma\right] + \frac{b}{\tau}\left[ -e_k \sin \sigma + f_k \cos \sigma\right] \nonumber \\
&& + (g + g^{sh}/2) e_m  - (b + b^{sh}/2) f_m, \quad \mbox{and} \\
\operatorname{Im} I_{mk} & = & -\frac{g}{\tau}\left[  -e_k \sin \sigma + f_k \cos \sigma\right] - \frac{b}{\tau} \left[e_k \cos \sigma + f_k \sin \sigma\right] \nonumber \\
&& + (b + b^{sh}/2) e_m + (g + g^{sh}/2) f_m.
\end{eqnarray}
Thus,
\begin{eqnarray}
P_{mk} & = & \operatorname{Re} V_m I^*_{mk} \nonumber \\
& = & e_m \left\{ -\frac{g}{\tau}\left[ e_k \cos \sigma + f_k \sin \sigma\right] + \frac{b}{\tau}\left[ -e_k \sin \sigma + f_k \cos \sigma\right] \right\} \nonumber \\
&& + f_m \left\{ -\frac{g}{\tau}\left[  -e_k \sin \sigma + f_k \cos \sigma\right] - \frac{b}{\tau} \left[e_k \cos \sigma + f_k \sin \sigma\right] \right\} \nonumber \\
&& + (g + g^{sh}/2) (e^2_m + f^2_m) \nonumber \\
& = & \left[ e_m - \frac{1}{\tau} e_k \cos \sigma \right](g \, , \, b) (^{e_m}_{f_m}) + \left[ f_m - \frac{1}{\tau} f_k \cos \sigma \right](-b \, , \, g) (^{e_m}_{f_m}) \nonumber \\
&& + \frac{g^{sh}}{2}(e_m^2 + f_m^2) \nonumber \\
&& + \frac{ -g e_m f_k - b f_m f_k - b e_m e_k + g f_m e_k }{\tau} \sin \sigma \nonumber \\
& = & \left[ e_m - \frac{1}{\tau} e_k \cos \sigma - \frac{1}{\tau} f_k \sin \sigma\right](g \, , \, b) (^{e_m}_{f_m}) + \left[ f_m - \frac{1}{\tau} f_k \cos \sigma + \frac{1}{\tau} e_k \sin \sigma\right](-b \, , \, g) (^{e_m}_{f_m}) \nonumber \\
&& + \frac{g^{sh}}{2}(e_m^2 + f_m^2). \label{rectpmk3}
\end{eqnarray}
\noindent We now turn to the representation of $P_{km}$ and 
$P_{mk}$ using polar coordinates.
\begin{eqnarray}
S_{km} & = & V_k I^*_{km} \ = \ V_k \left(\frac{1}{\tau^2}V_k^* - \frac{1}{\tau}e^{-j \sigma} V_m^* \right)y^* + \frac{1}{2 \tau^2} (g^{sh} -j b^{sh}) V_k V_k^* \label{skmpolar3a} \\
& = & \ V_k \left(\frac{1}{\tau^2}V_k^* - \frac{1}{\tau}e^{-j \sigma} V_m^* \right)(g - jb) + \frac{1}{2 \tau^2} (g^{sh} -j b^{sh}) V_k V_k^*\nonumber \\
& = & | V_k |^2 \frac{g}{\tau^2}  -  |V_k||V_m|\frac{g}{\tau} \cos (\theta_{km} - \sigma) -  |V_k||V_m|\frac{b}{\tau} \sin(\theta_{km} - \sigma)  + \frac{g^{sh}}{2\tau^2}|V_k|^2\nonumber \\
&& + j \left[ -|V_k|^2 \frac{b}{\tau^2} + |V_k| |V_m| \frac{b}{\tau} \cos(\theta_{km} - \sigma) - |V_k||V_m| \frac{g}{\tau} \sin(\theta_{km} - \sigma) - \frac{b^{sh}}{2\tau^2}|V_k|^2\right]. \nonumber
\end{eqnarray}
Similarly,
\begin{eqnarray}
S_{mk} & = & | V_m |^2 {g}  -  |V_k||V_m|\frac{g}{\tau} \cos (\theta_{mk} + \sigma) -  |V_k||V_m|\frac{b}{\tau} \sin(\theta_{mk} + \sigma)  + \frac{g^{sh}}{2}|V_m|^2\nonumber \\
&& + j \left[ -|V_m|^2 b + |V_k| |V_m| \frac{b}{\tau} \cos(\theta_{mk} + \sigma) - |V_k||V_m| \frac{g}{\tau} \sin(\theta_{mk} + \sigma) - \frac{b^{sh}}{2}|V_m|^2\right]. \nonumber
\end{eqnarray}
Hence the active power loss equals
\begin{eqnarray}
&& P_{km} + P_{mk} \ = \ \nonumber \\
&& \left(\frac{|V_k|^2}{\tau^2} + |V_m|^2\right) g - 2 g\frac {|V_k|}{\tau}|V_m| \cos(\theta_{km} - \sigma) +  \frac{g^{sh}}{2\tau^2}|V_k|^2 \ + \ \frac{g^{sh}}{2}|V_m|^2. \label{loss3a}
\end{eqnarray}
There is an alternative derivation of this equation that proves useful.  Consider point $k_1$ in Figure \ref{pimodel}.  The power injection into the line, at $k_1$, is equal to $P_{km}$ (i.e. it equals $\frac{V_k}{N} (N*I_{km})^* = V_k I_{km}^* = P_{km}$.  Moreover by construction the voltage magnitude at $k_1$ equals $|V_k|/\tau$ and the phase angle difference from $k_1$ to $m$ equals $\theta_{km} - \sigma$.  We can now recover (\ref{loss3a}) from (\ref{lossmain}), with the 
last two terms account for shunts, as when deriving (\ref{loss2}).

\subsubsection {$\Delta$ and loss inequalities, 3}
In the transformer case there will be two $\Delta$-inequalities. The first is obtained by from (\ref{rectpkm3}) 
by taking absolute values:
\begin{eqnarray}
&& |P_{km}| - \frac{g^{sh}}{2 \tau^2} |V_k|^{2} \ \le \ \nonumber \\
&& \frac{\mu_{km}}{\tau} \left|\frac{e_k}{\tau} - e_m \cos \sigma + f_m \sin \sigma  \right|  \, + \, \frac{\nu_{km}}{\tau} \left |\frac{1}{\tau} f_k - f_m \cos \sigma - e_m \sin \sigma  \right|. \label{delta3km}
\end{eqnarray}
Here as before $\mu_{km}$ and $\nu_{km}$ are known upper bounds on
$| (g \, , \,  b) (^{e_k}_{f_k})|$ and $| (-b \, , \, g) (^{e_k}_{f_k})|$, 
respectively.  Similarly, we obtain a  second $\Delta$-inequality 
from (\ref{rectpmk3}):
\begin{eqnarray}
&& |P_{mk}| - \frac{g^{sh}}{2 \tau^2} |V_m|^2 \ \le \ \nonumber \\
&& \mu_{mk} \left|e_m - \frac{1}{\tau} e_k \cos \sigma - \frac{1}{\tau} f_k \sin \sigma \right|  \, + \, \nu_{mk} \left | f_m - \frac{1}{\tau} f_k \cos \sigma + \frac{1}{\tau} e_k \sin \sigma \right|. \label{delta3mk}
\end{eqnarray}
Thus in order to represent these inequalities we need to introduce additional
lifted variables used to model $\left|\frac{e_k}{\tau} - e_m \cos \sigma + f_m \sin \sigma  \right|$, $\left|e_m - \frac{1}{\tau} e_k \cos \sigma - \frac{1}{\tau} f_k \sin \sigma \right|$, $\left |\frac{1}{\tau} f_k - f_m \cos \sigma + e_m \sin \sigma  \right|$ and $\left | f_m - \frac{1}{\tau} f_k \cos \sigma + \frac{1}{\tau} e_k \sin \sigma \right|$.  In the no-transformer case the first two
variables are equal to  $|e_m - e_k|$ and the last two are equal to $|f_m - f_k|$.  Replacing, in (\ref{delta3km}) and (\ref{delta3mk}), $|V_k|^2$ and $|V_m|^2$ with 
$V_k^{(2)}$ and $V_m^{(2)}$ respectively, we obtain the most general form of the $\Delta$-inequalities. \\

To obtain a loss inequality we apply the reasoning following equation (\ref{loss3a}).  Note that the voltage at point $k_1$ satisfies
\begin{eqnarray}
V_{k_1} & = & \frac{V_k}{N} \ = \ \frac{1}{\tau}(e_k + j f_k)(\cos \sigma - j \sin \sigma) \ = \ \nonumber \\
&& \frac{1}{\tau}( e_k \cos \sigma + f_k \sin \sigma) \, + \, j \frac{1}{\tau}(f_k \cos \sigma - e_k \sin \sigma). 
\end{eqnarray}
Since all losses are incurred in the section of the line between $k_1$ and $m$, 
applying (\ref{loss2}) we obtain:
\begin{eqnarray}
&& g \left(e_m  -  \frac{1}{\tau} \cos \sigma - \frac{1}{\tau} f_k \sin \sigma \right)^2 + g \left(f_m - \frac{1}{\tau} f_k \cos \sigma + \frac{1}{\tau} e_k \sin \sigma \right)^2  \nonumber \\ 
&& \ \le \ P_{km} + P_{mk} - \frac{g^{sh}}{2} \left(\frac{V_k^{(2)}}{\tau^2} + V_m^{(2)}\right). \label{loss3full1}
\end{eqnarray}
In this form we obtain a convex inequality that employs the auxiliary variables
introduced in (\ref{delta3mk}).  A similar construction yields an inequality
using the auxiliary variables in (\ref{delta3km}). 

\subsubsection {Circle inequalities, 3}
In the transformer case the structure of the circle inequalities differs due to
the asymmetry caused by the transformer.  First, the system (\ref{circlesystem2}) applied at $m$ is unchanged (i.e. system (\ref{circlesystem2}) with $k$ and 
$m$ interchanged).  To obtain a system at $k$ we again consider point $k_1$ in
Figure (\ref{pimodel}) and we now obtain:
\begin{subequations}
\label{circlesystem3}
\begin{eqnarray}
P_{km} - \left( \frac{r }{r^2 + x^2} + \frac{g^{sh}}{2} \right)  \frac{V_k^{(2)}}{\tau^2}  & = & \alpha_{km} \\
Q_{km} - \left( \frac{x }{r^2 + x^2}  - \frac{b^{sh}}{2}\right) \frac{V_k^{(2)} }{\tau^2} & = & \beta_{km} \\
\alpha_{km}^2 + \beta_{km}^2 & \le & \frac{V_k^{(2)} V_m^{(2)} }{\tau^2(r^2 + x^2)}.
\end{eqnarray}
\end{subequations}

\section{Inequalities from semidefinite relaxations}
Let $w$ be a subvector of the vector with entries $(1, e_1, e_2, \ldots, e_N, f_1, f_2, \ldots, f_N)^T$ where $N$ is the number of buses.  Then we can insist that $ww^T \succeq 0$. 
This is a semidefinite constraint.  As an alternative, by introducing additional lifted variables
we can use linear separating inequalities as a substitute for the semidefinite 
constraint.    Each additional lifted variable corresponds to an entry of the matrix $ww^T$.

To fix ideas, suppose that $w = (e_1, e_3, f_2)^T$.  Then we will have lifted variables
corresponding $e_1^2, e_3^2, f_2^2, e_1e_3, e_1f_2, e_3f_2$, which we denote, respectively,
by $e_1^s, e_3^s, f_2^s, ee_{1e}, ef_{12}$ and $ef_{32}$. Suppose that particular
values of these variables are such that the resulting matrix
\begin{eqnarray}
\tilde W & = & \left( \begin{array}{lcr}
\tilde e_1^s &  \tilde ee_{13} & \tilde ef_{12}\\
\tilde ee_{13}& \tilde e^s_{3}& \tilde ef_{32}\\
\tilde ef_{12} & \tilde ef_{32} & \tilde f_2^s
\end{array} \right) \ \not \succeq 0.
\end{eqnarray}
Let $u \in \R^3$ be such that $u^T \tilde W u < 0$.  [Given $\tilde W$ such a vector can
be computed in polynomial time, for example by running an adaptation of the Cholesky 
factorization procedure, or directly by computing an appropriate eigenvector of $\tilde W$.]  Then the inequality 
\begin{eqnarray}
u_1^2 \, e_1^s + 2 u_1 u_2 \, ee_{13} + 2 u_1u_3 \, ef_{12} + u_2^2 \, e^s_e + 2 u_2u_3 \, ef_{32} + u_3^2 \, f_s^s \ \ge \ 0
\end{eqnarray}
is valid, and violated by $\tilde W$.  We term such an inequality a \textit{semidefinite cut}.  In computation, the (perhaps obvious) requirement that $\norm{u} = 1$ is important,
as is the appropriate choice of bus indices used to construct the vector $w$.

\section{Tightening inequalities through reference angle fixings}\label{tighten}
Above we introduced a family of inequalities for each line of the underlying
network.  Here we will describe a tightening procedure that can render significant improvements.

Recall the discussion in Section \ref{basic} regarding foundational and lifted 
variables. Thus,
the lifted variables include e.g. a variable used to represent the
quantity $\left|e_m - \frac{1}{\tau} e_k \cos \sigma - \frac{1}{\tau} f_k \sin \sigma \right|$ introduced in equation (\ref{delta3mk}). 

We can express these facts in compact form as follows.  As in Section \ref{basic}, let $(P, Q, V^{(2)})$ indicate the vector of all foundational variables. Here, for each bus $k$ 
variable $V^{(2)}_k$ is
used to represent the quantity $|V_k|^2$.  If $N$ and $M$ indicate the 
number of buses and lines, respectively, then $(P, Q, V^2) \in \R^{2M + N}$.  
Let $W$ indicate the vector of all lifted variables, say with $H$ components, 
and let $K \subseteq \R^{2M + N + H}$ indicate
the convex set described by all inequalities introduced above.  Then we 
can represent $(P, Q, V^{(2)}, W) \in K$ more compactly by stating that
\begin{eqnarray}
 (P, Q, V^{(2)}) \, \in \, \hat K \, \doteq \, {\proj}_{\R^{2M + N}} K \label{projection}
\end{eqnarray}
where $\proj_{\R^{2M + N}} K$ is the projection of $K$ to the subspace of
the first $2M + N$ variables.

We now describe a procedure for tightening (\ref{projection}). As is well known, fixing an arbitrary bus at an arbitrary angle does not change
the set of feasible solutions to a standard OPF problem.  Thus, let $\hat k$ be
a particular bus, and let $\hat \theta_{\hat k}$ be a particular angle; we can therefore without loss of generality fix $\theta_{\hat k} = \hat \theta_{\hat k}$.  How can we take advantage of this fact so as to obtain stronger constraints?  Trivially, we can
of course enforce $f_{\hat k} = \tan\hat \theta_{\hat k} e_{\hat k}$.  

Moreover, consider for example the $\Delta$-inequality (\ref{delta1}) for a line
$\hat k m$ (for simplicity we assume the line has zero shunt admittance and no 
transformer).  We repeat the constraint here for convenience:
\begin{eqnarray}
|P_{\hat km}| & \le & \mu_{\hat km} | e_{\hat k} - e_m|  \, + \, \nu_{\hat km} | f_{\hat k} - f_m|, \label{delta1b}
\end{eqnarray}
where $\mu_{\hat km}$  and $\nu_{\hat km}$ are valid upper bounds on
$$ | (g \, , \,  b) (^{e_{\hat k}}_{f_{\hat k}})| \quad \mbox{and} \quad | (-b \, , \, g) (^{e_{\hat k}}_{f_{\hat k}})|,$$
respectively. [As previously both $b$ and $g$ depend on the line but we omit the dependency
for simplicity of notation].  Given that we know $\theta_{\hat k} = \hat \theta_{\hat k}$
we can tighten the estimates on $\mu_{\hat km}$  and $\nu_{\hat km}$, thereby obtaining
a tighter inequality from (\ref{delta1b}).  We can likewise tighten many of the 
inequalities introduced above.

More generally, suppose that rather than fixing $\theta_{\hat k}$ to a fixed value, we 
insist that it is contained in a known set $I(\hat k)$ (in particular an interval), i.e.
$$\theta_{\hat k} \, \in \, I^{(\hat k)}$$
As just argued we can therefore \textit{without loss of generality}, tighten the valid inequalities we described in previous section.  [This tightening is easiest in the case where
the set is in fact an interval.]  Let $K\left(\hat k, I^{(\hat k)}\right) \subseteq \R^{2M + N + H}$ denote the resulting convex body, and let 
$$ \Pi\left(\hat k, I^{(\hat k)}\right) \, \doteq \, {\proj}_{\R^{2M + N}} K\left(\hat k, I^{(\hat k)}\right).$$
As a consequence of the above observations, we now formally have:
\begin{LE} \label{fix1} Suppose $(\tilde P, \tilde Q, \tilde V^{(2)})$ is feasible for the OPF problem.  Then for any
bus $\hat k$, and any set $I^{(\hat k)}$,
\begin{eqnarray}
&& (\tilde P, \tilde Q, \tilde V^{(2)}) \ \in \ \Pi\left(\hat k, I^{(\hat k)}\right). \label{projk}
\end{eqnarray}
\qed
\end{LE}
Of course one can simply enforce (\ref{projk}) by explicitly writing down all the lifted variables and all the constraints used to describe the set $K\left(\hat k, I^{(\hat k)}\right)$.  Alternatively, one can \textit{separate} from the convex set $\Pi\left(\hat k, I^{(\hat k)}\right)$ and use such cuts as cutting planes.  From this perspective, the following
result is important:
\begin{CO} \label{fix1} Suppose $(\tilde P, \tilde Q, \tilde V^{(2)})$ is feasible for the OPF problem.  Then for any family of buses $k_i$ ($i \in F$) and sets $I^{(k_i)}$ we have
\begin{eqnarray}
&& (\tilde P, \tilde Q, \tilde V^{(2)}) \ \in \ \bigcap_{i \in F} \Pi\left(k_i, I^{(k_i)}\right). \label{projk}
\end{eqnarray}
\qed
\end{CO}
In other words, in particular, we can separate a given vector $(\tilde P, \tilde Q, \tilde V^{(2)})$ from sets obtained from our original family of valid inequalities by e.g. fixing
one arbitrary bus to an arbitrary angle, and tightening.

\section{Lower and upper bounds through linear mixed-integer programming techniques}\label{mip}
To address a the OPF problem in rectangular coordinates we use a technique that was 
originally developed by Glover \cite{glover}, used in \cite{mine} and more recently
analyzed in \cite{santanu}.  Suppose $u$ and $v$ are real variables.  We wish to approximate the product $uv$ with linear inequalities, to arbitrary precision.  Such a goal can
be achieved by adding a moderate number of binary variables.

To that effect, assume without loss of generality that $0 \le u \le 1$, $0 \le v \le 1$.  This can be achieved by translating and scaling the original $u$ and $v$, if they are assumed bounded.  Let $T \ge 1$ be an integer.  Then we can write
\begin{eqnarray}
 u  & = &  \sum_{j = 1}^T 2^{-j} y_j \ + \ \delta \label{expansion}
\end{eqnarray}
where each $y_j$ takes value zero or one, and $0 \le \delta \le 2^{-T}$.  Consequently,
we can approximate
\begin{eqnarray}
&& \sum_{j = 1}^T 2^{-j} y_j v \ \le \ uv \ \le \sum_{j = 1}^T 2^{-j} y_j v \ + \ 2^{-T} v \label{firstone}
\end{eqnarray}
This expression cannot directly be used because of the bilinear terms $y_j v$.  However,
let us set $w_j \doteq y_j v$.  Then we have
\begin{subequations}
\label{bilinear}
\begin{eqnarray}
w_j & \le & \min\{ v, y_j \} \\
w_j & \ge & \max\{v + y_j - 1, 0\}.
\end{eqnarray}
\end{subequations}
If $y_j = 1$ this system implies $w_j = v$ whereas if $y_j = 0$ the system yields $w_j = 0$. Hence the system, over binary $y_j$ but continuous $w_j$ and $v$ this system is a valid
relaxation of the bilinear relationship $w_j = y_j v$.   Thus, system (\ref{bilinear})
together with
\begin{eqnarray}
&& \sum_{j = 1}^T 2^{-j} w_j \ \le \ uv \ \le \sum_{j = 1}^T 2^{-j} w_j \ + \ 2^{-T} v \label{firstoneplus} \\
&& y_j \, \in \, \{0, 1\}, \quad 1 \le j \le K
\end{eqnarray}
yields an approximation to the quantity $uv$.  The bounds in (\ref{firstoneplus}) can 
be used appropriately to substitute each instance of the product $uv$ in an optimization
problem with a linear expression.

By performing the binary expansion
(\ref{expansion}) for selected rectangular coordinates $e_k$ or $f_k$ (suitably 
translated and rescaled) all bilinearities in the rectangular OPF formulation are removed.
The resulting optimization problem can then be run using a standalone mixed-integer 
solver; in particular with the goals of attaining good upper bounds
(albeit modulo the approximation errors
of magnitude $2^{-K}$), improving lower
bounds and possibly even proving infeasibility (same caveat as before regarding 
approximation errors).

\section{Initial computational experiments}
In the experiments reported here, we implemented the $Delta$, loss and circle inequalities
in their most general form.  To solver conic and linear programs, we used Gurobi v. 5.6.3
\cite{gurobi}.  To solve semidefinite programs, we used the system due to Lavaei and
coauthors \cite{javad2}, which also includes a procedure for extracting a feasible 
rank-one solution from the SDP.  All runs were performed on a current workstation with ample physical memory.  All running times are in seconds unless indicated.

In the table ``SDP time'' is the time taken to solve the SDP relaxation of the OPF problem,
``SDP gap'' is the percentage gap between the value of the SDP relaxation and the upper
bound (value of feasible solution) obtained by the SDP system.  ``SOCP time'' and ``LP time'', are, respectively, the time required to solve our conic relaxation and its first-order
(outer) relaxation through a cutting-plane algorithm.  ``SOCP gap'' and ``LP gap'' are
the percentage gaps relative to the SDP upper bound.

\begin{table}[ht]
\begin{center}
\begin{tabular}{|l|c|c|c|c|c|c|} \hline
 & {\bf SDP time} & {\bf SDP gap} & {\bf SOCP time} & {\bf SOCP gap} & {\bf LP time} & {\bf LP gap} \\ \hline
{\bf case9} & 1.04 & 0.0002 $\%$& 0.05 & 0.7899 $\%$ & 0.04 & 0.7899 $\%$ \\
{\bf case30} & 3.40 & 0.0185 $\%$& 0.23 & 1.3808 $\%$ & 0.35 & 1.3964 $\%$ \\
{\bf case57} &  4.23 & 0.0000 $\%$&  0.62 & 0.9954 $\%$ & 1.41 & 0.9954 $\%$ \\
{\bf case118} & 8.73 & 0.0045 $\%$& 0.98 & 1.4645 $\%$ &5.12  & 1.4642 $\%$ \\
{\bf case300} & 20.29  & 0.0018 $\%$& 4.62 & 1.0585 $\%$ & 49.61 & 1.0559 $\%$ \\
{\bf case2383wp} & 13 min & 0.6836 $\%$& 2 min & 3.6134 $\%$ & 1.63 & 5.6489  $\%$ \\
{\bf case2746wp} & 16 min & 0.0375 $\%$& 79.10 & 1.8593 $\%$ & 1.88 & 3.1235 $\%$ \\
\hline
\end{tabular}
\caption{Initial computational results.}
\label{firsttable}
\end{center}
\end{table}

\noindent \tiny    Tue.Nov..4.192606.2014@littleboy
\end{document}